# Spanning trees short or small[*]

R. Ravi[†]   R. Sundaram[‡]   M. V. Marathe[§]   D. J. Rosenkrantz[¶]

S. S. Ravi[‖]

**Abstract**

We study the problem of finding small trees. Classical network design problems are considered with the additional constraint that only a specified number $k$ of nodes are required to be connected in the solution. A prototypical example is the $k$MST problem in which we require a tree of minimum weight spanning at least $k$ nodes in an edge-weighted graph. We show that the $k$MST problem is NP-hard even for points in the Euclidean plane. We provide approximation algorithms with performance ratio $2\sqrt{k}$ for the general edge-weighted case and $O(k^{1/4})$ for the case of points in the plane. Polynomial-time exact solutions are also presented for the class of decomposable graphs which includes trees, series-parallel graphs, and bounded bandwidth graphs, and for points on the boundary of a convex region in the Euclidean plane.

We also investigate the problem of finding short trees, and more generally, that of finding networks with minimum diameter. A simple technique is used to provide a polynomial-time solution for finding $k$-trees of minimum diameter. We identify easy and hard problems arising in finding short networks using a framework due to T. C. Hu.

[*]This paper appeared in a preliminary form as [28].

[†]Dept. of Computer Science, University at California at Davis, CA 95616. Email: ravi@cs.ucdavis.edu. Research done while author was at Brown University with the support of an IBM Graduate Fellowship, NSF PYI award CCR-9157620 and DARPA contract N00014-91-J-4052 ARPA Order No. 8225.

[‡]Dept. of Computer Science, MIT LCS, Cambridge MA 02139. Email: koods@theory.lcs.mit.edu. Research supported by DARPA contract N0014-92-J-1799 and NSF 92-12184 CCR.

[§]Dept. of Computer Science, SUNY, Albany, NY 12222. Email: madhav@cs.albany.edu. Supported by NSF Grant CCR 89-03319.

[¶]Dept. of Computer Science, SUNY, Albany, NY 12222. Email: djr@cs.albany.edu. Supported by NSF Grant CCR 90-06396.

[‖]Dept. of Computer Science, SUNY, Albany, NY 12222. Email: ravi@cs.albany.edu. Supported by NSF Grant CCR 89-05296.

# 1 Introduction

## 1.1 Motivation: small trees

The oil reconnaissance boats are back from their final trip off the coast of Norway[1], and present you with a detailed map of the seas surrounding the coastline. Marked in this map are locations which are believed to have a good chance of containing oil under the sea bed. Your company has a limited number of oil rigs that it is willing to invest in the effort. Your problem is to position these oil rigs at marked places so that the cost of laying down pipelines between these rigs is minimized. The problem at hand can be modeled as follows: Given an edge-weighted graph and a specified number $k$, find a tree of minimum weight spanning at least $k$ nodes. Note that a solution to the problem will be a tree spanning exactly $k$ nodes. We call this problem the *k-Minimum Spanning Tree* (or the $k$MST) problem. In this paper, we study such classical network-design problems as the MST problem with the additional constraint that only a specified number of nodes need to be incorporated into the network. Unlike the MST problem which admits a polynomial-time solution [4, 22, 25], the $k$MST problem is considerably harder to solve[2].

**Theorem 1.1** *The $k$MST problem is NP-complete*[3].

The above theorem holds even when all the edge weights are drawn from the set $\{1, 2, 3\}$ (or any set containing three distinct values). It is not hard to show a polynomial-time solution for the case of two distinct weights. The problem remains NP-hard even for the class of planar graphs as well as for points in the plane.

## 1.2 Approximation algorithms

A $\rho$-approximation algorithm for a minimization problem is one that delivers a solution of value at most $\rho$ times the minimum. Consider a generalization of the $k$MST problem, the $k$-Steiner tree problem: given an edge-weighted graph, an integer $k$ and a subset of at least $k$ vertices specified as terminals, find a minimum-weight tree spanning at least $k$ terminals. We can apply approximation results for the $k$MST problem to this problem by considering the auxiliary complete graph on the terminals with edges weighted by shortest-path distances. A $\rho$-approximation for the $k$MST problem on the auxiliary graph yields a $2\rho$-approximation for the $k$-Steiner tree problem. Therefore we focus on approximations for the $k$MST problem. We provide the first approximation algorithm for this problem.

**Theorem 1.2** *There is a polynomial-time algorithm that, given an undirected graph $G$ on $n$ nodes with nonnegative weights on its edges, and a positive integer $k \leq n$, constructs a tree spanning at least $k$ nodes of weight at most $2\sqrt{k}$ times that of a minimum-weight tree spanning any $k$ nodes.*

The algorithm in the above theorem is based on a combination of a greedy technique that constructs trees using edges of small cost and a shortest-path heuristic that merges trees when the number of trees to be merged is small. The analysis of the performance ratio is based on a solution-decomposition technique [10, 21, 26, 27] which uses the structure of the optimal solution to derive a bound on the cost of the solution constructed by the approximation algorithm.

---
[1]Story reconstructed from a communication from Naveen Garg [16].
[2]The main theorems in this paper are stated in the introduction and proved in later sections.
[3]This result was independently obtained by Lozovanu and Zelikovsky [23].

The above theorem provides a $4\sqrt{k}$-approximation algorithm for the $k$-Steiner tree problem as well. Moreover, we can construct an example that demonstrates that the performance guarantee of the approximation algorithm is tight to within a constant factor.

We can derive a better approximation algorithm for the case of points in the Euclidean plane.

**Theorem 1.3** *There is a polynomial-time algorithm that, given $n$ points in the Euclidean plane, and a positive integer $k \leq n$, constructs a tree spanning at least $k$ of these points such that the total length of the tree is at most $O(k^{\frac{1}{4}})$ times that of a minimum-length tree spanning any $k$ of the points.*

As before, we can construct an example showing that the performance ratio of the algorithm in Theorem 1.3 is tight. Our proof of Theorem 1.3 also yields as corollary an approximation algorithm for the rectilinear $k$MST problem.

**Corollary 1.4** *There is a polynomial-time algorithm that, given $n$ points in the plane, and a positive integer $k \leq n$, constructs a rectilinear tree spanning at least $k$ of these points such that the total length of the tree is at most $O(k^{\frac{1}{4}})$ times that of a minimum-length rectilinear tree spanning any $k$ of the points.*

## 1.3 Exact algorithms: special cases

Since the $k$MST problem is NP-complete even for the class of planar graphs, we focus on special classes of graphs and provide exact solutions that run in polynomial time. Bern, Lawler and Wong [7] introduced the notion of decomposable graphs. A class of decomposable graphs is defined using a finite number of primitive graphs and a finite collection of binary composition rules. Examples of decomposable graphs include trees[4], series-parallel graphs and bounded-bandwidth graphs. We use a dynamic programming technique to prove the following theorem.

**Theorem 1.5** *For any class of decomposable graphs, there is an $O(nk^2)$-time algorithm for solving the $k$MST problem.*

Though the $k$MST problem is hard for arbitrary configurations of points in the plane, we have the following result.

**Theorem 1.6** *There is a polynomial-time algorithm for solving the $k$MST problem for the case of points in the Euclidean plane that lie on the boundary of a convex region.*

---

[4]A polynomial-time algorithm for trees was also independently obtained by Lozovanu and Zelikovsky [23].

The proof of the above theorem uses a monotonicity property of the optimal tree along with a degree constraint on an optimal solution. This allows us to apply dynamic programming to find the exact solution. Several researchers in computational geometry have presented exact algorithms for choosing $k$ points that minimize other objectives such as diameter, perimeter, area and volume [2, 12, 13, 14].

## 1.4 Short trees

Keeping the longest path in a network small is often an important consideration in network design. We investigate the problem of finding networks with small diameter. Recall that the diameter of a tree is the maximum distance (path length) between any pair of nodes in the tree. The problem of finding a minimum-diameter spanning tree of an edge-weighted graph was shown to be polynomially solvable by Camerini, Galbiati and Maffioli [9] when the edge weights are nonnegative. They also show that the problem becomes NP-hard when negative weights are allowed. Camerini and Galbiati [8] have proposed polynomial-time algorithms for a bounded path tree problem on graphs with nonnegative edge weights. Their result can be used to show that the minimum-diameter spanning tree problem as well as its natural generalization to Steiner trees can be solved in polynomial time. We use a similar technique to show that the following *minimum-diameter $k$-tree* problem is polynomially solvable: given a graph with nonnegative edge weights, find a tree of minimum diameter spanning at least $k$ nodes.

**Theorem 1.7** *There is a polynomial-time algorithm for the minimum-diameter $k$-tree problem on graphs with nonnegative edge weights.*

We investigate easy and hard results in finding short networks. For this, we use a framework due to T. C. Hu [19]. In this framework, we are given a graph with nonnegative distance values $d_{ij}$ and nonnegative requirement values $r_{ij}$ between every pair of nodes $i$ and $j$ in the graph. The communication cost of a spanning tree is defined to be the sum over all pairs of nodes $i, j$ of the product of the distance between $i$ and $j$ in the tree under $d$ and the requirement $r_{ij}$. The objective is to find a spanning tree with minimum-communication cost. Hu considered the case when all the $d$ values are one and showed that a Gomory-Hu cut tree [18] using the $r$ values as capacities is an optimal solution. Hu also considered the case when all the $r$ values are one and derived sufficient conditions under which the optimal tree is a star. The general version of the latter problem is NP-hard [9, 20].

We define the diameter cost of a spanning tree to be the maximum cost over all pairs of nodes $i, j$ of the distance between $i$ and $j$ in the tree multiplied by $r_{ij}$. In Table 1, we present current results in this framework. All $r_{ij}$ and $d_{ij}$ values are assumed to be nonnegative. The first two rows of the table examine the cases when either of the two parameters is uniform-valued. The last two rows illustrate that the two problems become NP-complete when both the parameters are two-valued.

## 1.5 Short small trees

We consider the $k$-tree versions of the minimum-communication-cost and minimum-diameter-cost spanning tree problems and show the following hardness result.

| $r_{ij}$ | $d_{ij}$ | Communication cost | Diameter cost |
|---|---|---|---|
| Arbitrary | $\{a\}$ | Cut-tree [19] | Open |
| $\{a\}$ | Arbitrary | NP-complete [20] | Poly-time [9] |
| $\{a,b\}$ | $\{0,c\}$ | Cut-tree variant (this paper) | Poly-time (this paper) |
| $\{a,4a\}$ | $\{c,d\}$ | NP-complete [20] | NP-complete (this paper) |

Table 1: Results on minimum-communication-cost spanning trees and minimum-diameter-cost spanning trees.

**Theorem 1.8** *The minimum-communication $k$-tree problem and the minimum-diameter $k$-tree problem are both hard to approximate within any factor even when all the $d_{ij}$ values are one and the $r_{ij}$ values are nonnegative.*

In the next section, we present the NP-completeness results. Section 3 contains the $2\sqrt{k}$ approximation for the $k$MST problem. In Section 4, we present the stronger result for the case of points in the plane. Then we address polynomially solvable cases of the problem. In Section 6, we prove our results on short trees. We close with a discussion of directions for future research.

## 2 NP-completeness results

In this section we show that the $k$MST problem is NP-hard by reducing the Steiner tree problem to it. The Steiner tree problem is known to be NP-hard [15]. As an instance of the Steiner tree problem we are given an undirected graph $G$, a set of terminals $R$ (which is a subset of the vertex set of $G$) and a positive integer $M$, and the question is whether there exists a tree spanning $R$ and containing at most $M$ edges. We transform this input to an instance $G', k$, of the $k$MST problem as follows: We let $X = |V(G)| - |R| + 1$ and connect each terminal of $G$ to a distinct path of $X$ new vertices, the path consisting of zero-weighted edges. We assign weight one to the already existing edges of $G$ and set the weight between all other pairs of vertices to $\infty$ (a very large number). This is the graph $G'$ (See Figure 1). We set $k$ to be $|R| \cdot X$. If there exists a Steiner tree in $G$ spanning the set $R$ and containing at most $M$ edges, then it is easy to construct a $k$MST of weight at most $M$ in $G'$. Conversely, by our choice of $k$ and $X$, any $k$MST in $G'$ must contain at least one node from the path corresponding to each terminal in $R$. Hence any $k$MST can be used to derive a Steiner tree for $R$ in $G$. This completes the reduction. Extensions of hardness to the case of planar graphs and points in the plane follow in a similar way from the hardness of the Steiner tree problem in these restricted cases. Given a planar embedding of $G$ we can create an embedded version of $G'$ since only paths are added.

The NP-hardness holds even when all the edge costs are from the set $\{1,2,3\}$. The reduction for this case is similar to the above. Without loss of generality we assume that in the given instance of the Steiner tree problem, $G$ is connected and $M \leq |V| - 1$. We let $X = |V(G)| - |R| + 1$ as before, and connect each terminal of $G$ to a distinct set of $X$ vertices by edges of weight one. We set the original edges of $G$ to have weight two and all other edges to have weight three. We choose $k = |R| \cdot X + M + 1$ and the bound on the cost of the $k$MST to be $|R| \cdot X + 2M$. If there exists a Steiner tree in $G$ spanning the set $R$ and containing at most $M$ edges, then it is easy to construct a $k$MST of weight at most $|R| \cdot X + 2M$ in $G'$. This is done by connecting all the newly added vertices to the Steiner tree using the weight one edges and then picking up more vertices (note that the graph is connected and $M \leq |V| - 1$) using the weight two edges until there are $|R| \cdot X + M + 1$ vertices. If there exists a $k$MST of weight at most $|R| \cdot X + 2M$ in $G'$ then note that the $k$MST cannot contain an edge of weight three because it has only $k - 1 = |R| \cdot X + M$ edges and if it contained an edge of weight three then it would have to contain at least $|R| \cdot X + 1$ edges of weight one but there are only $|R| \cdot X$ edges of weight one in $G'$. Further, the $k$MST must span $R$, and since it has at most $M$ edges of weight two, hence there must exist a Steiner tree in $G$ spanning $R$ and containing at most $M$ edges.

When there are only two distinct edge costs, the $k$MST problem can be solved in polynomial time. The basic idea is the following: Let $w_1$ and $w_2$ denote the two edge weights, where $w_1 < w_2$. Construct an edge subgraph $G_1$ of $G$ containing all the edges of weight $w_1$. Choose a minimum number, say $r$, of the connected components of $G_1$ to obtain a total of $k$ nodes. Construct a spanning tree for each chosen component and connect the trees together into a single tree by adding exactly $r - 1$ edges of weight $w_2$. It is straightforward to verify that the resulting solution is optimal.

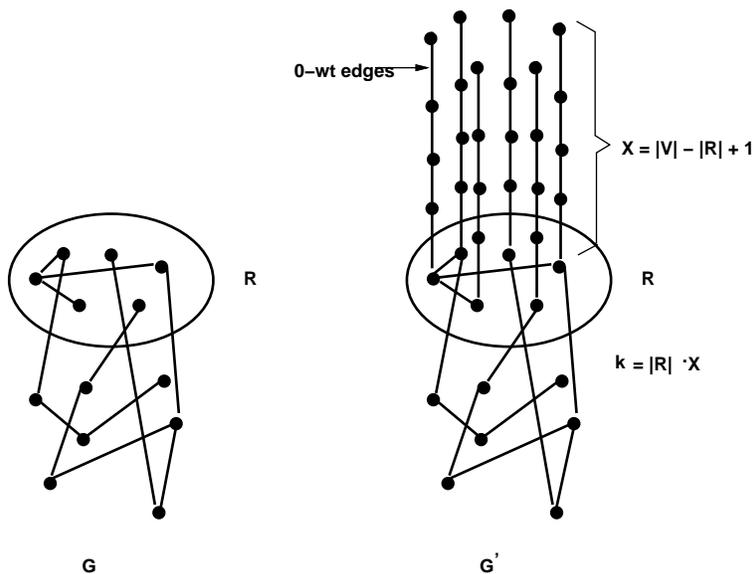

Figure 1: The basic NP-hardness reduction from Steiner tree to $k$MST.

# 3 The approximation algorithm for the general case

In this section, we present the proof of Theorem 1.2. As input, we are given an undirected graph $G$ with nonnegative edge weights and an integer $k$.

## 3.1 The algorithm and its running time

It is useful to think of the algorithm as running in two distinct phases: a merge phase and a collect phase.

During the merge phase, the algorithm maintains a set of clusters and a spanning tree on the vertex set of each cluster. Initially each vertex forms a singleton cluster. At each step of the merge phase, we choose an edge of minimum cost among all edges that are between two clusters, and merge them by using the edge to connect their spanning trees.

Define the size of a cluster to be the number of vertices that it contains. During the course of the merge phase, the clusters grow in size. The collect phase is entered only when

(i) there exist at most $\sqrt{k}$ clusters whose sizes sum to at least $k$, and

(ii) no cluster has size $k$ or more.

In the collect phase, we consider each cluster in turn as the root and perform a shortest-path computation between clusters using the weights on inter-cluster edges. We determine for each cluster $C$, the shortest distance $d_C$ such that, within distance $d_C$ from $C$, there exist at most $\sqrt{k}$ clusters whose sizes sum to at least $k$. Note that by the first precondition for starting the collect phase, the distance $d_C$ is well defined. We choose the cluster $C$ with the minimum value of $d_C$ and connect it using shortest paths of length at most $d_C$ to each of these $\sqrt{k}$ clusters. We can prune edges from some of these shortest paths to output a tree of clusters whose sizes sum to $k$. We may do this since any cluster has less than $k$ nodes at the start of this phase by the second precondition.

The merge phase of the algorithm continues to run until both the preconditions of the collect phase are satisfied. Beginning with the step of the merge phase after which both preconditions of the collect phase are satisfied, at each subsequent step, the algorithm forks off an execution of the collect phase for the current configuration of clusters. The merge phase continues to run until a cluster of size $k$ or more is formed. Next, merge phase prunes the edges of the spanning tree of the cluster whose size is between $k$ and $2k$ so as to obtain a spanning tree of size exactly $k$. At this point, the merge phase terminates and outputs the spanning tree of the cluster of size $k$. Each forked execution of the collect phase outputs a spanning tree of size between $k$ and $2k$ as well. The algorithm finally outputs the tree of least weight among all these trees. The algorithm is given below:

**Algorithm Merge-Collect**

1. Initialize each vertex to be in singleton connected components and the set of edges chosen by the algorithm to be $\phi$. Initialize the iteration count $i = 1$.

2. Repeat until there exists a cluster whose size is between $k$ and $2k$

(a) Let $VS_i = \{C_1 \cdots C_l\}$ denote the set of connected components at the start of this iteration. Assume that the components are numbered in non-increasing order of their size.

(b) Form an auxiliary graph $G(VS_i, E')$ where the edge $(C_i, C_j)$ between two components is the minimum cost edge in $E$ whose endpoints belong to $C_i$ and $C_j$ respectively.

(c) Choose a minimum cost edge $(C_i, C_j)$ in $G(VS_i, E')$ and merge the corresponding clusters $C_i$ and $C_j$.

(d) $VS_{i+1} = VS_i - \{C_i\} - \{C_j\} \cup \{C_i \cup C_j\}$

**Remark:** This corresponds to one iteration of merge phase.

(e) Let $j^* = \min\{j : \sum_{i=1}^{j} |C_i| \geq k\}$.

(f) If $j^* \leq \sqrt{k}$ then $SOL_i = \text{Collect}(G(VS, E'))$

(g) $i = i + 1$;

3. Prune the edges of the cluster whose size is between $k$ and $2k$ to obtain a tree with exactly $k$ vertices. Denote the tree obtained by $MSOL$.

4. The output of the heuristic is the minimum valued tree among $MSOL$ and all the $SOL_i$'s.

**Procedure Collect$(G(V, E))$**

1. For each cluster vertex $C$ do

   (a) With the cluster $C$ as the root, form a shortest path tree.

   (b) Let $d_C$ denote the shortest distance from $C$ such that there are no more than $\sqrt{k}$ clusters whose sizes sum up to at least $k$.

   (c) Choose these clusters and join them to the root cluster by using the edges in the shortest path tree computed in Step 1(a).

   (d) Prune the edges of the tree to obtain a tree having exactly $k$ nodes.

2. Output the tree corresponding to the choice of the root cluster $C$ that minimizes $d_C$.

It is easy to see that there are at most $O(n)$ steps in the merge phase and hence at most this many instances of the collect phase to be run. Using Djikstra's algorithm [11] in each collect phase, the whole algorithm runs in time $O(n^2(m + n \log n))$ where $m$ and $n$ denote the number of edges and nodes in the input graph respectively. The running time of the collect phase dominates the running time of the merge phase.

## 3.2 The performance guarantee

Consider an optimal $k$MST of weight $OPT$. During the merge phase, nodes of this tree may merge with other nodes in clusters. We focus our attention on the number of edges of the optimal $k$MST that are *exposed*, i.e., remain as inter-cluster edges. We show that at any step in which a large number of edges of the $k$MST are exposed, every edge in the spanning tree of each cluster has small weight.

**Lemma 3.1** *If at the beginning of a step of the merge phase, an optimal $k$MST has at least $x$ exposed edges (inter-cluster edges), then each edge in the spanning tree of any cluster at the end of the step has weight at most $\frac{OPT}{x}$.*

**Proof:** The proof uses induction on the number of steps. Suppose that an optimal $k$MST has at least $x$ exposed edges at the beginning of the current step of the merge phase. Then at the beginning of the previous step, the optimal $k$MST must have had at least $x$ exposed edges as well. Thus by the induction hypothesis every edge in the spanning tree of any cluster at the end of the previous step has weight at most $\frac{OPT}{x}$. Since only one new composite cluster is formed in the current step, it remains to show that the edge added in this iteration has cost at most $\frac{OPT}{x}$. But this is straightforward since there is an optimal $k$MST with at least $x$ exposed edges of total weight at most $OPT$. □

We now prove the performance guarantee in Theorem 1.2. The above lemma is useful as long as the number of exposed edges is high. Applying the lemma with $x = \sqrt{k}$ shows that every edge in the spanning tree of each cluster has weight at most $\frac{OPT}{\sqrt{k}}$. Consider the scenario when the merge phase runs to completion to produce a tree with at least $k$ nodes even before the number of exposed edges falls below $\sqrt{k}$. In this case, since the resulting tree has at most $k$ nodes, the cost of the tree is at most $\frac{OPT}{\sqrt{k}} \cdot k \leq 2\sqrt{k} \cdot OPT$.

Otherwise, the number of exposed edges falls below $\sqrt{k}$ before the merge phase runs to completion. However, in this case, note that both preconditions for the start of the collect phase will have been satisfied. Hence the algorithm must have forked off a run of the collect phase. We show that the tree output by this run has low weight. Consider a shortest-path computation of the collect phase rooted at a cluster containing a node of the optimal $k$MST. Then clearly, within a distance at most $OPT$, we can find at most $\sqrt{k}$ clusters whose sizes sum to at least $k$. Since the number of exposed edges is less than $\sqrt{k}$, the clusters containing nodes of the optimal tree form such a collection. Since there are at most $\sqrt{k}$ clusters to connect to, the weight of these connections is at most $\sqrt{k} \cdot OPT$. It remains to bound the weight of the spanning trees within each of the clusters retained in the output solution. This is not hard since all edges in these clusters have weight at most $\frac{OPT}{\sqrt{k}}$ by Lemma 3.1. Since the size of the output tree is at most $k$ (as a result of the pruning), the total weight of all the edges retained within these clusters is at most $\sqrt{k} \cdot OPT$. Summing the weight of these intra-cluster edges and the inter-cluster connections shows that the output tree has cost at most $2\sqrt{k} \cdot OPT$. This proves the performance ratio of $2\sqrt{k}$ claimed in Theorem 1.2.

The example in Figure 2 shows that the performance ratio of the algorithm is $\Omega(\sqrt{k})$.

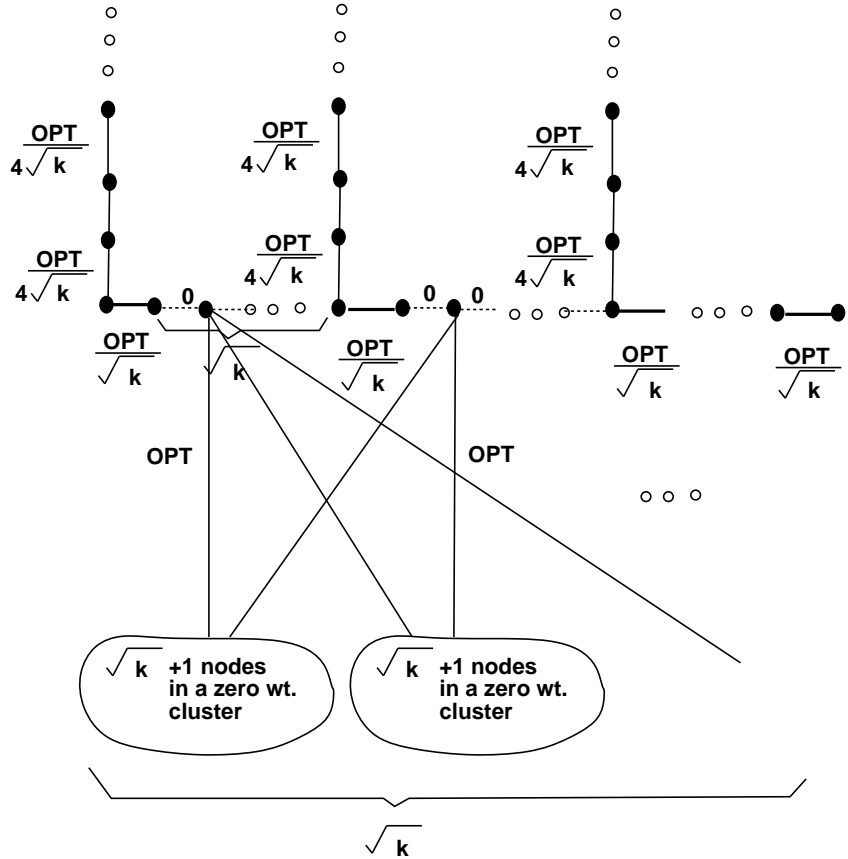

Figure 2: Example of a graph in which the algorithm in Theorem 1.2 outputs a tree of weight $\Omega(OPT \cdot \sqrt{k})$. The optimal $k$MST is the horizontal path made of zero-weight edges and the $\sqrt{k}$ edges of weight $\frac{OPT}{\sqrt{k}}$ each. All zero-weight edges will be chosen first in the merge phase. The merge phase running to completion will extend each of the zero-weight upward-directed paths to include $\Omega(k)$ edges each of weight $\frac{OPT}{4\sqrt{k}}$ resulting in a tree of weight $\Omega(OPT \cdot \sqrt{k})$. The collect phases may output trees consisting of all the $\sqrt{k}+1$-sized clusters at the bottom of the figure each of weight $\Omega(OPT \cdot \sqrt{k})$.

# 4 An approximation algorithm for points on the plane

In this section, we present a heuristic for the $k$MST problem for points on the plane and a proof of its performance guarantee. Let $S = \{s_1, s_2, ..., s_n\}$ denote the given set of points. For any pair of points $s_i$ and $s_j$, let $d(i, j)$ denote the Euclidean distance between $s_i$ and $s_j$.

## 4.1 The heuristic

I. **For** each distinct pair of points $s_i$, $s_j$ in $S$ **do**

  (1) Construct the circle $C$ with diameter $\delta = \sqrt{3}d(i,j)$ centered at the midpoint of the line segment $\langle s_i, s_j \rangle$.

  (2) Let $S_C$ be the subset of $S$ contained in $C$. If $S_C$ contains fewer than $k$ points, skip to the next iteration of the loop (i.e., try the next pair of points). Otherwise, do the following.

  (3) Let $Q$ be the square of side $\delta$ circumscribing $C$.

  (4) Divide $Q$ into $k$ square cells each with side $= \delta/\sqrt{k}$.

  (5) Sort the cells by the number of points from $S_C$ they contain and choose the minimum number of cells so that the chosen cells together contain at least $k$ points. If necessary, arbitrarily discard points from the last chosen cell so that the total number of points in all the cells is equal to $k$.

  (6) Construct a minimum spanning tree for the $k$ chosen points. (For the rectilinear case, construct a rectilinear minimum spanning tree for the $k$ chosen points.)

  (7) The solution value for the pair $\langle s_i, s_j \rangle$ is the length of this MST.

II. **Output** the smallest solution value found.

It is easy to see that the above heuristic runs in polynomial time. In the next subsection, we show that the heuristic provides a performance guarantee of $O(k^{1/4})$. We begin with some lemmas.

## 4.2 The performance guarantee

**Lemma 4.1** *Let $S$ denote a set of points on the plane, with diameter $\Delta$. Let $a$ and $b$ be two points in $S$ such that $d(a,b) = \Delta$. Then the circle with diameter $\sqrt{3}\Delta$ centered at the midpoint of the line segment $\langle a, b \rangle$ contains $S$.*

**Proof:** Suppose there exists a point $p \in S$ not contained within the circle of diameter $\sqrt{3}\Delta$ centered at the midpoint of the line segment $\langle a, b \rangle$. If $p$ lies on the perpendicular bisector of the line segment $\langle a, b \rangle$ then it is clear that $d(a,p) = d(b,p) > \Delta$, else $p$ is closer to one of $a$ and $b$ than the other. Say $p$ is closer to $a$; then it is easy to see that $d(b,p) > \Delta$. Thus, if there exists a point outside the circle then it contradicts the fact that the diameter of the set $S$ is $\Delta$. Hence $S$ must be contained within the circle. □

**Lower Bounds on an Optimal $k$MST**

The following lemma is used to establish a lower bound on $OPT$.

**Lemma 4.2** *Consider a square grid on the plane with the side of each cell being $\sigma$. Then the length of an MST for any set of $t$ points, where each point is from a distinct cell is $\Omega(t\sigma)$.*

**Proof:** Pick a point from the set and discard all points in the eight cells neighboring the cell containing the chosen point. Doing this repeatedly we choose a subcollection of $t/9 = \Omega(t)$ points such that the distance between any pair of points in the subcollection is at least $\sigma$. The lemma then follows from the observation that the minimum length of a tree spanning $\Omega(t)$ points that are pairwise $\sigma$-distant is $\Omega(t\sigma)$. □

Let $P^*$ denote the set of points in an optimal solution to the problem instance. Let $\Delta$ denote the *diameter* of $P^*$ (i.e., the maximum distance between a pair of points in $P^*$), and $OPT$ denote the length of an MST for $P^*$. Consider an iteration in which the circle constructed by the heuristic is defined by two points $a$ and $b$ in $P^*$ such that $d(a,b) = \Delta$. Let $g$ be the number of square cells used by the heuristic in selecting $k$ points in this iteration. To establish the performance guarantee of the heuristic, we show that the length of the MST constructed by the heuristic during this iteration is within a factor $O(k^{1/4})$ of $OPT$.

It is easy to see that $OPT \geq \Delta$ because $\Delta$ is the diameter of $P^*$.

Since the heuristic uses a minimum number ($g$) of square cells in selecting $k$ points, the points in $P^*$ must occur in $g$ or more square cells. Note that the side of each square cell is $\sqrt{3}\Delta/\sqrt{k}$. This gives us the following corollary to Lemma 4.2.

**Corollary 4.3**
$$OPT = \Omega(g\Delta/\sqrt{k})$$

**Upper Bound on the Cost of the Heuristic**

We now prove an upper bound on the cost of the spanning tree returned by the heuristic.

For this, we need the following lemma.

**Lemma 4.4** *The length of a minimum spanning tree for any set of $q$ points in a square with side $\sigma$ is length $O(\sigma\sqrt{q})$.*

**Proof:** Paste a square grid over the square where each sub-cell in the grid has side $\sigma/\sqrt{q}$. Connect each point to a closest vertex in the grid. Consider the tree consisting of one vertical line, all the horizontal lines in the grid connected to the vertical line, and the vertical lines connecting each point to its nearest horizontal line (See Figure 3). It is clear that the grid lines in the tree have total length $O(\sigma\sqrt{q})$ and the lines connecting the points to the grid have total length $q \cdot O(\sigma/\sqrt{q}) = O(\sigma\sqrt{q})$. □

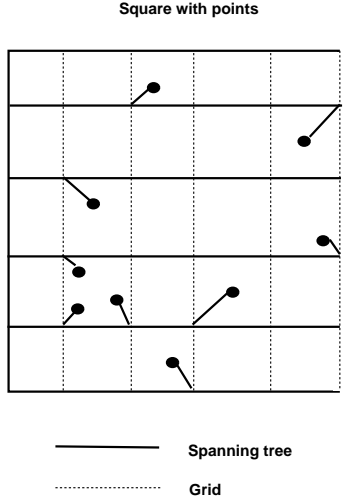

Figure 3: A spanning tree of length $O(\sigma\sqrt{q})$ on any $q$ points in a square of side $\sigma$.

**Lemma 4.5** *The length of the spanning tree constructed by the heuristic is $O(\sqrt{g}\Delta)$.*

**Proof:** Let $Q_i$ denote the set of points in the $i^{th}$ cell chosen by the heuristic, $1 \leq i \leq g$. Thus $\sum_{i=1}^{g} |Q_i| = p$. Consider the following two-stage procedure for constructing a spanning tree for the points in $\cup_{i=1}^{g} Q_i$.

*Stage I:* Construct a minimum spanning tree for the points in $Q_i$, $1 \leq i \leq g$. Note that the points in $Q_i$ are within a square of side $\sqrt{3}\Delta/\sqrt{p}$. Using Lemma 4.4, the length of an MST for $Q_i$ is $O(\frac{\Delta}{\sqrt{k}}\sqrt{|Q_i|})$. Thus, the total length of all the minimum spanning trees constructed in this stage is $O(\frac{\Delta}{\sqrt{p}}\sum_{i=1}^{g}\sqrt{|Q_i|}) = O(\sqrt{g}\,\Delta)$ by the Cauchy-Schwartz inequality.

*Stage II:* Connect the $g$ spanning trees constructed in Stage I into a single spanning tree as follows. Choose a point arbitrarily from each $Q_i$ ($1 \leq i \leq g$), and construct an MST for the $g$ chosen points. Note that these $g$ points are within a square of side $\sqrt{3}\,\Delta$. Thus, by Lemma 4.4, the length of the MST constructed in this stage is $O(\sqrt{g}\,\Delta)$ as well.

Thus, the total length of the spanning tree constructed by the two-stage procedure is $O(\sqrt{g}\,\Delta)$. □

**The Final Analysis**
We are now ready to complete the proof of the performance bound. As argued above, $OPT = \Omega(\Delta)$, and from Corollary 4.3, $OPT = \Omega(g\Delta/\sqrt{k})$. Thus $OPT = \Omega(\max\{\Delta, g\Delta/\sqrt{k}\})$. Also from Lemma 4.5, the length of the spanning tree produced by the heuristic is $O(\sqrt{g}\,\Delta)$. Therefore, the performance ratio is $O(\min\{\sqrt{g}, \sqrt{k/g}\}) = O(k^{1/4})$ as claimed.

The example in Figure 4 shows that the performance ratio of the heuristic is $\Omega(k^{1/4})$.

Observe that both our lower bounds on an optimal solution and the upper bound on the spanning tree obtained also apply to the case of constructing a rectilinear $k$MST. Hence it follows that the above approximation algorithm delivers a performance guarantee of $O(k^{1/4})$ for the rectilinear $k$MST problem too. This proves Theorem 1.4.

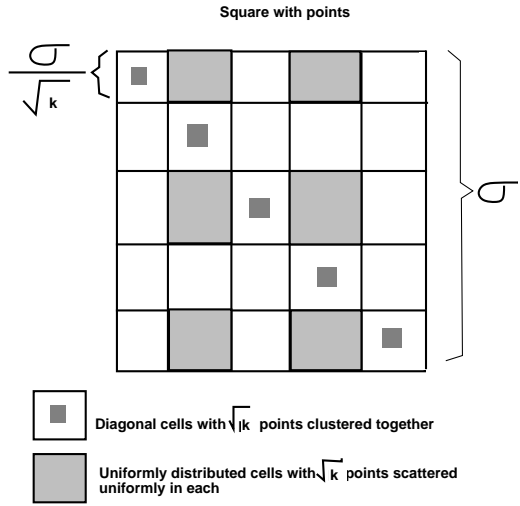

Figure 4: Example of a configuration of points on the plane in which the heuristic outputs a tree of length $\Omega(OPT \cdot \sqrt{k})$. The big square has side $\sigma$. Each cell of the square grid has side $\sigma/\sqrt{k}$. There are $\sqrt{k}$ points clustered closely together in each cell along the diagonal of the big square. And in each of $\sqrt{k}$ cells distributed uniformly throughout the big square there are $\sqrt{k}$ uniformly distributed points. The heuristic may pick up the points in the uniformly distributed cells forming a tree of length $\Omega(\sigma \cdot k^{1/4})$ while the tree spanning the points along the diagonal has length $O(\sigma)$.

## 5 Exact algorithms for special cases

### 5.1 $k$MST for Decomposable Graphs

In this section, we prove Theorem 1.5. A class of **decomposable graphs** $\Gamma$ is given by a set of rules satisfying the following conditions [7].

1. The number of primitive graphs in $\Gamma$ is finite.

2. Each graph in $\Gamma$ has an ordered set of special nodes called **terminals**. The number of terminals in each graph is bounded by a constant.

3. There is a finite collection of binary composition rules that operate only at terminals, either by identifying two terminals or adding an edge between terminals. A composition rule also determines the terminals of the resulting graph, which must be a subset of the terminals of the two graphs being composed.

Examples of decomposable graphs include trees, series-parallel graphs, bounded-bandwidth graphs, etc. [7].

Let $\Gamma$ be any class of decomposable graphs. The $k$MST problem for $\Gamma$ can be solved optimally in polynomial time using dynamic programming. Following [7], it is assumed that a given graph $G$ is accompanied by a parse tree specifying how $G$ is constructed using the rules and that the size of the parse tree is linear in the number of nodes of $G$.

Consider a fixed class of decomposable graphs $\Gamma$. Suppose that $G$ is a graph in $\Gamma$. Let $\pi$ be a partition of a nonempty subset of the terminals of $G$. We define the following set of costs for $G$.

$Cost_i^\pi(G)$ = Minimum total cost of any forest containing a tree for each block of $\pi$, such that the terminal nodes occurring in each tree are exactly the members of the corresponding block of $\pi$, no pair of trees is connected, the total number of edges in the forest is $i$ and each tree contains at least one edge ($1 \leq i < k$).

$Cost_{k-1}^\emptyset(G)$ = Minimum cost of a tree within $G$ containing $k-1$ edges, and containing *no* terminal nodes of $G$.

For any of the above costs, if there is no forest satisfying the required conditions, the value of $Cost$ is defined to be $\infty$.

Note that because $\Gamma$ is fixed, the number of cost values associated with any graph in the parse tree for $G$ is $O(k)$. We now show how the cost values can be computed in a bottom-up manner, given the parse tree for $G$.

To begin with, since $\Gamma$ is fixed, the number of primitive graphs is finite. For a primitive graph, each cost value can be computed in constant time, since the number of forests to be examined is fixed. Now consider computing the cost values for a graph $G$ constructed from subgraphs $G_1$ and $G_2$, where the cost values for $G_1$ and $G_2$ have already been computed.

Let $\Pi_{G_1}$, $\Pi_{G_2}$ and $\Pi_G$ be the set of partitions of a subset of the terminals of $G_1$, $G_2$ and $G$ respectively. Let $A$ be the set of edges added to $G_1$ and $G_2$ by the composition rule $R$ used in constructing $G$ from $G_1$ and $G_2$. Corresponding to rule $R$, there is a partial function $f_R : \Pi_{G_1} \times \Pi_{G_2} \times 2^A \to \Pi_G$, such that a forest corresponding to partition $\pi_1$ in $\Pi_{G_1}$, a forest corresponding to partition $\pi_2$ in $\Pi_{G_2}$, and a subset $B \subseteq A$, combine to form a forest corresponding to partition $f_R(\pi_1, \pi_2, B)$ of $G$. Furthermore, if the forest corresponding to $\pi_1$ contains $i$ edges, and the forest corresponding to $\pi_2$ contains $j$ edges, then the combined forest in $G$ contains $i + j + |B|$ edges.

Similarly, there is a partial function $g_R : \Pi_{G_1} \times 2^A \to \Pi_G$, such that a forest corresponding to partition $\pi_1$ in $\Pi_{G_1}$ and a subset $B \subseteq A$ combine to form a forest corresponding to partition $g_R(\pi_1, B)$ of $G$. If the forest corresponding to $\pi_1$ contains $i$ edges, then the combined forest in $G$ contains $i + |B|$ edges. There is also a similar partial function $h_R : \Pi_{G_2} \times 2^A \to \Pi_G$. Finally, there is a partial function $j_R : 2^A \to \Pi_G$.

Using functions $f_R$, $g_R$, $h_R$ and $j_R$, cost values for $G$ can be computed from the set of cost values for $G_1$ and $G_2$. For instance, suppose that $f_R(\pi_1, \pi_2, B) = \pi$. Then a contributor to computing $Cost_i^\pi(G)$ is $Cost_t^{\pi_1}(G_1) + Cost_{i-t-|B|}^{\pi_2}(G_2) + w(B)$, for each $t$ such that $1 \leq t \leq i - |B| - 1$. Here $w(B)$ is the total cost of all edges in $B$. The value of

$Cost_i^\pi(G)$ is the minimum value among its contributors.

When all the cost values for the entire graph $G$ have been computed, the cost of an optimal $k$MST is equal to $\min_{\pi \in \pi_G}\{Cost_{k-1}^\pi(G)\}$, where the forest corresponding to $\pi$ consists of a single tree.

We now analyze the running time of the algorithm. For each graph occurring in the parse tree, there are $O(k)$ cost values to be computed. Each of the cost values can be computed in $O(k)$ time. As in [7], we assume that the size of the given parse tree for $G$ is $O(n)$. Then the dynamic programming algorithm takes time $O(nk^2)$. This completes the proof of Theorem 1.5.

### 5.2 $k$MST for points on the boundary of a convex region

We now restrict our attention to the case where we are given $n$ points that lie on the boundary of a convex region, and show that the $k$MST on these points can be computed in polynomial time using dynamic programming. We also provide a faster algorithm if the points are constrained to lie on the boundary of a circle.

**Lemma 5.1** *Any optimal $k$MST for a set of points in the plane is non self-intersecting.*

**Proof:** Suppose an optimal $k$MST were self intersecting, then let $\langle a,b \rangle$ and $\langle c,d \rangle$ be the intersecting line segments. On removing the edges $\langle a,b \rangle$ and $\langle c,d \rangle$ from the $k$MST we get three connected components, hence some two vertices, one from $\{a,b\}$ and one from $\{c,d\}$ must be in the same connected component. Say, $a$ and $d$ are in the same connected component, then since in any convex quadrilateral the sum of two opposite sides is less than the sum of the two diagonals, replacing $\langle a,b \rangle$ and $\langle c,d \rangle$ by $\langle a,c \rangle$ and $\langle b,d \rangle$ we still get a tree spanning $k$ nodes but with lesser weight. This contradicts the fact that the $k$MST we started out with was optimal. Hence any optimal $k$MST on a set of points in the plane must be non self-intersecting. □

**Lemma 5.2** *Given $n$ points on the boundary of a convex polygon no vertex in an optimal $k$MST of these points has degree greater than 4.*

**Proof:** Suppose there is a vertex $v$ in an optimal $k$MST with degree greater than 4. Let $v_1, v_2, \ldots, v_d, d \geq 5$ be its neighbors in the optimal $k$MST as shown in the figure (See Fig. 5.). Using the well known fact that any convex polygon lies entirely on one side of a supporting line, we have that $\angle v_1vv_d \leq 180°$. By the pigeon-hole principle, there is an $i$ such that $\angle v_ivv_{i+1} \leq 180°/(d-1) < 60°, 1 \leq i \leq d-1$ since $d$ is at least 5. Thus in $\triangle v_ivv_{i+1}$, $\angle v_ivv_{i+1}$ is not the largest angle, and $v_iv_{i+1}$ is not the largest side. Therefore replacing the larger of $vv_i$ and $vv_{i+1}$ in the optimal $k$MST with $v_iv_{i+1}$ we obtain a tree with lesser weight, contradicting the assumption that the $k$MST was optimal. This completes the proof. □

We now characterize the structure of an optimal solution in the following decomposition lemma and use it to define the subproblems which we need to solve recursively using dynamic programming.

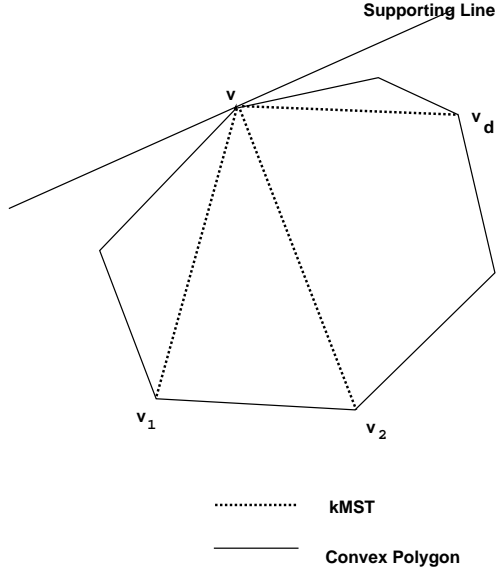

Figure 5: Points on a convex polygon.

**Lemma 5.3 (Decomposition lemma.)** *Let $v_0, v_1, \ldots, v_{n-1}$ be the vertices of a convex polygon in say, clockwise order. Let $v_i$ be a vertex of degree $d_i$ in an optimal kMST. Note that $1 \leq d_i \leq 4$.*

*If $d_i \geq 2$ let the removal of $v_i$ from the optimal kMST produce connected components $C_1, C_2, \ldots, C_{d_i}$ (See Fig 6.). Let $|C_i|$ denote the number of vertices in component $C_i$. Then there exists a partition of $v_{i+1}, v_{i+2}, \ldots, v_{i-1}$, (indices taken $\bmod\ n$), into $d_i$ contiguous subsegments $S_1, S_2, \ldots, S_{d_i}$ such that $\forall j, 1 \leq j \leq d_i$, the optimal kMST induced on $S_j \bigcup \{v_i\}$ is an optimal $(|C_j|+1)MST$ on $S_j \bigcup \{v_i\}$ in which the degree of $v_i$ is one.*

*If $d_i = 1$, let $v_j$ be $v_i$'s neighbor in the optimal kMST. Let $v_j$ be adjacent to $d_{j1}$ vertices in $v_{i+1}, v_{i+2} \ldots, v_{j-1}$ and $d_{j2}$ vertices in $v_{j+1}, v_{j+2}, \ldots, v_{i-1}$. Let the optimal kMST contain $|C_1|$ vertices from the set $v_{i+1}, v_{i+2} \ldots, v_{j-1}$ and $|C_2|$ vertices from the set $v_{j+1}, v_{j+2}, \ldots, v_{i-1}$. Then the optimal kMST induced on $v_{i+1}, v_{i+2} \ldots, v_j$ is an optimal $(|C_1|+1)MST$ on $v_{i+1}, v_{i+2} \ldots, v_j$ with degree of $v_j = d_{j1}$ and the optimal kMST induced on $v_j, v_{j+1} \ldots, v_{i-1}$ is an optimal $(|C_2|+1)MST$ on $v_j, v_{j+1} \ldots, v_{i-1}$ with degree of $v_j = d_{j2}$.*

**Proof:** If $d_i \geq 2$ then it is easy to see that a partition of $v_{i+1}, v_{i+2}, \ldots, v_{i-1}$ into contiguous subsegments $S_1, S_2, \ldots, S_{d_i}$ exists such that $\forall j, 1 \leq j \leq d_i, C_j \subset S_j$, because the optimal kMST is non self-intersecting by Lemma 5.1. Further, the optimal kMST induced on $S_j \bigcup \{v_i\}$ must be an optimal $(|C_j|+1)MST$ on $S_j \bigcup \{v_i\}$ with degree of $v_i = 1$, for otherwise we could replace it getting a lighter kMST. The proof of the case when $d_i = 1$ is equally straightforward and is omitted. □

Thus the subproblems we consider are specified by the following four parameters: a size $s$, a vertex $v_i$, the degree $d_i$ of $v_i$, and a contiguous subsegment $v_{k1}, v_{k1+1}, \ldots, v_{k2}$ such that $i \notin [k1 \ldots k2]$. A solution to such a subproblem denoted by $SOLN(s; v_i; d_i; v_{k1}, v_{k1+1}, \ldots, v_{k2})$ is the weight of an optimal $sMST$ on $\{v_i, v_{k1}, v_{k1+1}, \ldots, v_{k2}\}$ in which $v_i$ has degree $d_i$. Using the decomposition lemma above, we can write a simple recurrence relation for $SOLN(s; v_i; d_i; v_{k1}, v_{k1+1}, \ldots, v_{k2})$.

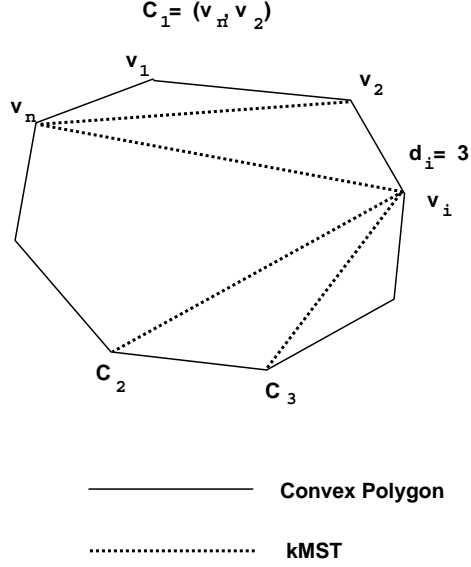

Figure 6: Decomposition.

$SOLN(s; v_i; d_i; v_{k1}, v_{k1+1}, \ldots, v_{k2}) =$

$$\begin{cases} \infty : \text{if } d_i = 0 \text{ or } s < d_i + 1 \text{ or } ((k2 - k1 + 1) \bmod n) + 1 < s. \\ \min_{k'_0 = k1 < k'_1 \ldots < k'_{d_i} = k2} \min_{s_1 \ldots + s_{d_i} = s + d_i - 1, s_j \geq 1} \Sigma_{1 \leq j \leq d_i} SOLN(s_j; v_i; 1; v_{k'_{j-1}}, \ldots, v_{k'_j}) : \text{if } d_i \geq 2 \\ \min_{j_0 = k1 \leq j_1 \leq j_2 = k2} \{w(v_i v_{j_1}) + \min_{0 \leq d_1 + d_2 \leq 3} \min_{s_1 + s_2 = s} \\ (SOLN(s_1; v_{j_1}; d_1; v_{j_0}, \ldots, v_{j_1 - 1}) + SOLN(s_2; v_{j_1}; d_2; v_{j_1 + 1}, \ldots, v_{j_2}))\}) : \text{if } d_i = 1 \end{cases}$$

Here $w(v_i v_j)$ is the cost of the edge $(v_i, v_j)$. The optimal $k$MST $=$

$$\min_{1 \leq i \leq n} \min_{1 \leq d \leq 4} SOLN(k; v_i; d; v_{i+1}, v_{i+2}, \ldots, v_{i-1})$$

Note that we have $O(kn^3)$ subproblems and each subproblem requires looking up the solution to at most $O(k^3 n^3)$ smaller subproblems. This yields a running time of $O(k^4 n^6)$. When $k = \Omega(\sqrt{n})$, this running time can be further improved by organizing the computation of the recurrences for the smaller subproblems better. Consider a vertex $v$, an integer $0 \leq s \leq k$ denoting the size of tree, and a partition of the other $(n-1)$ vertices into four groups. This corresponds to one of the subproblems we need to solve. Each smaller subproblem of this subproblem is specified by the number of nodes $s(\leq k)$ in the tree, vertex $v_i$ of degree $d_i$ in the interval $v_{k1}, \ldots, v_{k2}$ can be solved by first computing a partition of the interval into at most four parts (exactly four when $d_i = 4$). For the first subinterval, we compute the best tree containing $v$ on $i$ nodes with other nodes only from this subinterval, and $v$ has degree one in this tree, for $1 \leq i \leq s$. This computation takes $O(nk)$ times since there are at most $s \leq k$ trees to be computed, and for each $i$, there are at most $n$ nodes with which $v_i$ shares the single edge in the best tree. Next, we include the next subinterval, and compute for $1 \leq i \leq s$, the best tree on $i$ nodes containing $v_i$ and nodes from these two subintervals, where $v_i$ has degree two with one edge to a node in the first and one edge to a node in the second subinterval. This set of trees can also be computed in $O(nk)$ time given the set of trees for the first subinterval as follows: First, compute the best tree on $i$ nodes for $1 \leq i \leq s$ containing $v_i$ and nodes only in the second subinterval, where $v_i$ has exactly one edge to a node in this subinterval, in $O(nk)$ time as before. Using these values and the analogous set of values for the first subinterval, the best $i$ trees for the first two subintervals can be obtained in $O(k^2) = O(nk)$ time, since each of the $s \leq k$ trees requires looking up at most $s$ diffferent pairs of trees, one from each subinterval. This method can be extended to compute the solution for the whole set of four subintervals in $O(nk)$ time. Since there are $O(n^3)$ ways to partition a given interval into four subintervals, the recurrence for this subproblem can be solved in $O(kn^4)$ time. So the total time to solve one subproblem is $O(kn^4)$ time. Since there are a total of $O(kn^3)$ subproblems, the total running time of the algorithm is $O(k^2 n^7)$.

We now provide a faster algorithm to find the optimal $k$MST in the case when all $n$ points lie on a circle. We assume that no two points are diametrically opposite.

**Lemma 5.4** *Given $n$ points $v_1, v_2, \ldots, v_n$ on a circle no vertex in an optimal $k$MST has degree more than 2.*

**Proof:** Suppose point $v_p$ in an optimal $k$MST has degree greater than 2. Then consider the diameter passing through $v_p$. At least two neighbors of $v_p$ lie on one side of this diameter. Let these neighbors be $v_q$ and $v_r$, where $v_q$ is closer to $v_p$ than $v_r$. Then since $\angle v_p v_q v_r$ is obtuse we can replace $v_p v_r$ by $v_q v_r$ to get a smaller tree. □

Lemma 5.4 implies that if the points lie on a circle then every optimal $k$MST is a path. Moreover, if the path "zig-zags", then we can replace the crossing edge with a smaller edge. Thus we have the following lemma.

**Lemma 5.5** *Given $n$ points $v_1, v_2, \ldots, v_n$ on a circle, let a minimum length $k$-path on these points be $v_{i_1}, \ldots, v_{i_p}$. Then the line segment joining $v_{i_1}$ and $v_{i_p}$ along with the $k$-path forms a convex $k$-gon.*

**Proof:** By Lemma 5.4 the minimum-length $k$-path is also the minimum-length $k$MST. Suppose the line segment joining $v_{i_1}$ and $v_{i_p}$ along with the minimum $k$-path does not form a convex $k$-gon. Then there exists a zig-zag in the path as shown in Figure 7. Say the center of the circle lies to the right of the edge $\langle a, b \rangle$ then we can replace $\langle a, b \rangle$ by the edge $\langle b, c \rangle$ to get a smaller $k$MST which contradicts the fact that the $k$-path we started out with was optimal. □

Lemmata 5.4 and 5.5 lead to a straightforward dynamic programming algorithm to compute an optimal $k$MST for points on a circle in $O(n^3)$ time.

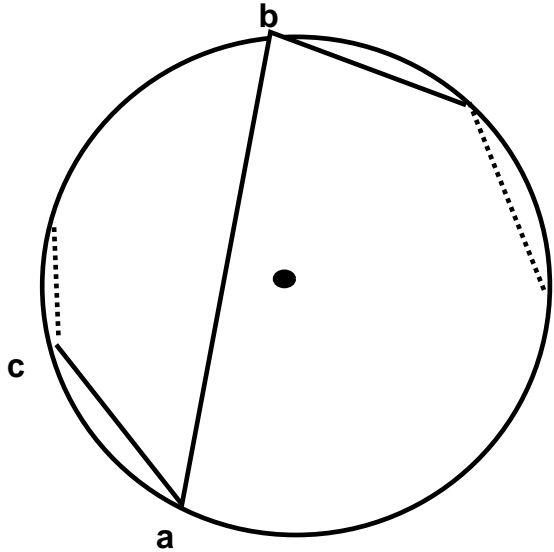

Figure 7: Illustration of Lemma 5.5.

# 6 Short trees and short small trees

## 6.1 Short trees

In this subsection, we prove our results on short trees. First, we address the minimum-diameter $k$-tree problem: Given a graph with nonnegative edge weights, find a tree of minimum diameter spanning at least $k$ nodes.

Recall that the diameter of a tree is the maximum distance (path length) between any pair of nodes in the tree. We introduce the notion of subdividing an edge in a weighted graph. A subdivision of an edge $e = (u, v)$ of weight $w_e$ is the replacement of $e$ by two edges $e_1 = (u, r)$ and $e_2 = (r, v)$ where $r$ is a new node. The weights of $e_1$ and $e_2$ sum to $w_e$. Consider a minimum-diameter $k$-tree. Let $x$ and $y$ be the endpoints of a longest path in the tree. The weight of this path, $D$, is the diameter of the tree. Consider the midpoint of this path between $x$ and $y$. If it falls in an edge, we can subdivide the edge by adding a new vertex as specified above. The key observation is that there exist at least $k$ vertices at a distance at most $D/2$ from this midpoint. This immediately motivates an algorithm for the case when the weights of all edges are integral and bounded by a polynomial in the number of nodes. In this case, all such potential midpoints lie in half-integral points along edges of which there are only a polynomial number. Corresponding to each candidate point, there is a smallest distance from this point within which there are at least $k$ nodes. We choose the point with the least such distance and output the breadth-first search (bfs) tree rooted at this point appropriately truncated to contain only $k$ nodes.

When the edge weights are arbitrary, the number of candidate midpoints are too many to check in this fashion. However, we can use a graphical representation of the distance of any node from any point along a given edge to bound the search for candidate points. We can think of an edge $e = (u, v)$ of weight $w$ as a straight line between its endpoints of length $w$. For any node $x$ in the graph, consider the shortest path from $x$ to a point along the edge $e$ at distance $\ell$ ($\leq w$) from $u$. The length of this path is the minimum of $\ell + d(x, u)$ and $w - \ell + d(v, x)$. We can plot this distance of the node $x$ as a function of $\ell$. The resulting plot is a piecewise linear bitonic curve that we call the roof curve of $x$ in $e$ (See Figure 8). For each edge $e$, we plot the roof curves of all the vertices of the graph in $e$. For any candidate point in $e$, the minimum diameter of a $k$-tree centered at this point can be determined by projecting a ray upwards from this point in the plot and determining the least distance at which it intersects the roof curves of at least $k$ distinct nodes. The best candidate point for a given edge is one with the minimum such distance. Such a point can be determined by a simple line sweep algorithm on the plot. Determining the best midpoint over all edges gives the midpoint of the minimum-diameter $k$-tree. This proves Theorem 1.7.

The following lemma gives yet another way to implement the polynomial time algorithm for finding a tree of minimum diameter spanning $k$ nodes.

**Lemma 6.1** *Given two vertices in a graph, $v_i$ and $v_j$, such that every other vertex is within distance $d_i$ of $v_i$ or $d_j$ of $v_j$, it is possible to find two trees, one rooted at $v_i$ and of depth at most $d_i$ and one rooted at $v_j$ of depth at most $d_j$ which partition the set of all vertices.*

**Proof:** Consider the shortest-path trees $T_i$ and $T_j$ rooted at $v_i$ and $v_j$ of depth $d_i$ and $d_j$, respectively. Every vertex occurs in one tree or both trees. Consider a vertex $v_p$ that occurs in both the trees. If it is the case that $d_i - \text{depth}_{T_i}(v_p)$ is greater than $d_j - \text{depth}_{T_j}(v_p)$ then the same is true of all descendants of $v_p$ in $T_j$. Hence we can remove $v_p$ and all it's descendants from $T_j$ since we are guaranteed that all these vertices occur in $T_i$. Repeating this procedure bottom-up we get two trees satisfying the required conditions and partitioning the vertex set. □

The above lemma motivates the following alternate algorithm for finding a minimum-diameter tree spanning at least $k$ nodes. For each vertex $v_i$ in the graph compute the shortest distance $d_i$ such that there are $k$ vertices within distance $d_i$ of $v_i$. For each edge $(v_i, v_j)$ compute the least $d_{ij}^i + d_{ij}^j$ such that there are $k$ vertices within distance $d_{ij}^i$ of $v_i$ or $d_{ij}^j$ of $v_j$. Then compute the least of all the $d_i$'s and $d_{ij}^i + d_{ij}^j + w(v_i, v_j)$'s and this is the diameter of the $k$-tree with least diameter.

We now address the results in the third row of Table 1.

**Lemma 6.2** *If the $r_{ij}$ values are drawn from the set $\{a, b\}$ and the $d_{ij}$ values from $\{0, c\}$ then the minimum-communication-cost spanning tree can be computed in polynomial time.*

**Proof:** When the $d_{ij}$ values are all uniform, Hu [19] observed that the Gomory-Hu cut tree with the $r_{ij}$ values as capacities is a minimum-communication-cost tree. We can use this result to handle the case when zero-cost $d_{ij}$ edges are allowed as well. We contract the connected components of the graph using zero-cost $d_{ij}$ edges into supernodes. The requirement value $r_{IJ}$ between two supernodes $v_I$ and $v_J$ is the sum of the requirement values $r_{ij}$ such that $i \in v_I$ and $j \in v_J$. Now we find a Gomory-Hu cut tree between the supernodes using the $r_{IJ}$ values as capacities. By choosing an arbitrary spanning tree of zero-$d_{ij}$-valued edges within each supernode and connecting them to the Gomory-Hu tree, we get a spanning tree of the whole graph. It is easy to verify that this is a minimum-communication-cost spanning tree in this case. □

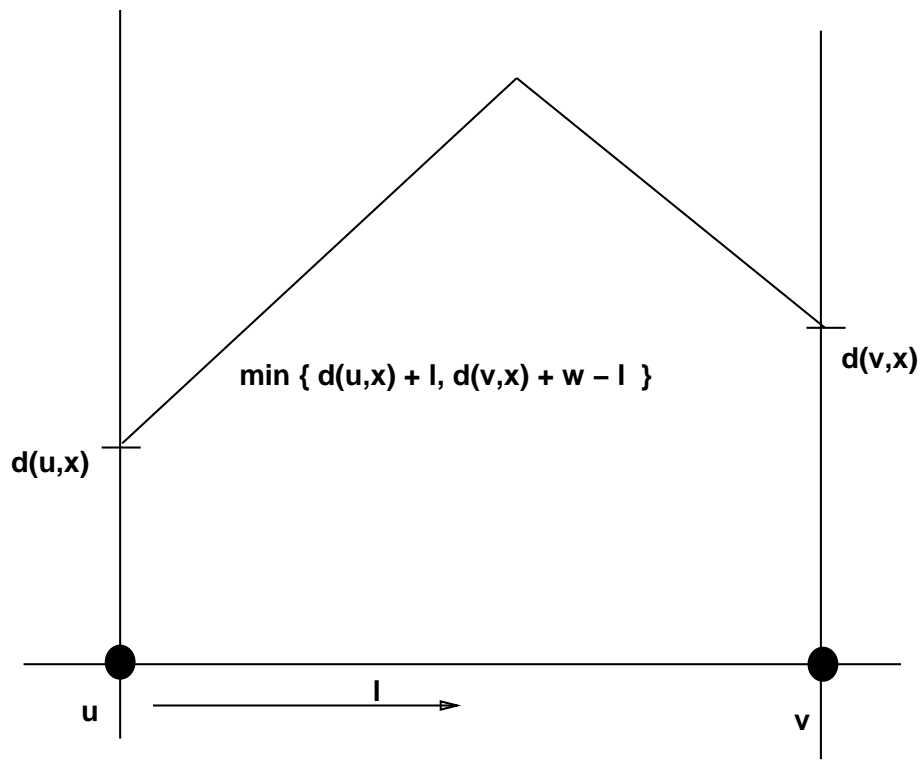

Figure 8: A roof curve of a node $x$ in edge $e = (u, v)$.

**Lemma 6.3** *When all the $d_{ij}$ values are uniform and there are at most two distinct $r_{ij}$ values (say a and b) then the minimum-diameter-cost spanning tree can be computed in polynomial time.*

**Proof:** Let the higher of the two $r_{ij}$ values be $a$. If the edges with requirement $a$ form a cyclic subgraph, then any spanning tree has diameter cost $2a$. In this case, any star is an optimal solution. Otherwise, consider the forest of edges with requirement $a$. Determine a center for each tree in this forest. Consider the tree formed by connecting these centers in a star. The root of the star is a center of the tree of largest diameter in the forest. If the diameter cost of the resulting tree is less than $2a$, it is easy to see that this tree has optimum diameter cost. Otherwise any star tree on all the nodes has diameter cost $2a$ and is optimal. Note that we can extend this solution to allow zero-cost $d_{ij}$ edges by using contractions as before. □

Now we address the results in the fourth row of Table 1.

**Lemma 6.4** *The minimum-diameter-cost spanning tree problem is NP-complete even when the $r_{ij}$'s and $d_{ij}$'s take on at most two distinct values.*

**Proof:** We use a reduction from an instance of 3SAT. We form a graph that contains a special node $t$ (the "true" node), a node for each literal and each clause. We use two $d_{ij}$ values, $c$ and $d$ where we assume $c < d$. Each literal is connected to its negation with an edge of distance $c$. The true node is connected to every literal with an edge of distance $c$. Each clause is connected to the three literals that it contains with edges of distance $c$. All other edges in the graph have distance $d$. Now we specify the requirements on the edges. We use requirement values from $\{a, 4a\}$, where $a \neq 0$. The requirement value of an edge between a literal and its negation is $4a$. The requirement value of all other edges is $a$ (See Figure 9). Assuming that $d > 4ac$, it is easy to check that there is a spanning tree of this graph with diameter cost at most $4ac$ if and only if the 3SAT formula is satisfiable. □

## 6.2 Short small trees

Finally we prove Theorem 1.8. We prove the theorem for the communication tree case. The proof of the other part is similar. Suppose there is a polynomial-time $M$-approximation algorithm for the minimum-communication-cost$k$-tree problem where all the $d_{ij}$ values are one and all $r_{ij}$ values are nonnegative. Then, we show that the $k$-independent set problem can be solved in polynomial time. The latter problem is well known to be NP-complete [15]. Given graph $G$ of the $k$-independent set problem, produce the following instance of the communication $k$-tree problem: $d_{ij} = 1$ for every pair of nodes $i, j$; Assign $r_{ij}$ equals one if $(i, j)$ is *not* an edge in $G$, and $Mk(k-1) + 1$ otherwise. If $G$ has an independent set of size $k$, then we can form a star on these $k$ nodes (choosing an arbitrary node as the root). In the star, the distance between any pair of nodes is at most 2 and the $r$ value for each pair is 1. Thus, the communication cost of an optimum solution is at most $k(k-1)$.

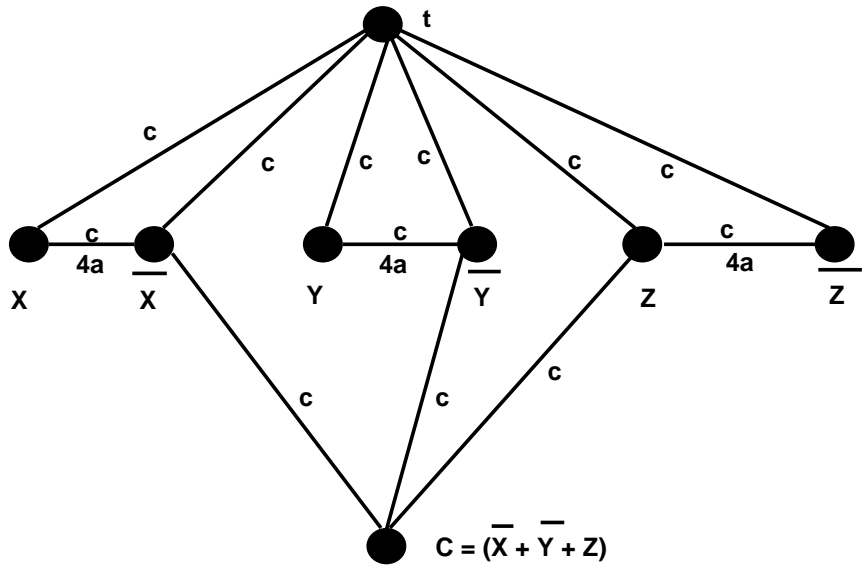

Figure 9: Reduction from an instance of 3SAT to the minimum-diameter-cost spanning tree problem.

The approximation algorithm will return a solution of cost at most $Mk(k-1)$. The nodes in this solution are independent in $G$ by the choice of $r_{ij}$ for nonedges $(i,j) \in G$. On the other hand, if there is no independent set of size $k$ in $G$, the communication cost of any $k$-tree is greater than $Mk(k-1)$.

# 7  Closing remarks

## 7.1  Future research

A natural question is whether there are approximation algorithms for the $k$MST problem which provide better performance guarantees than those presented in this paper. An interesting observation in this regard is the following. Any edge in an optimal $k$MST is a shortest path between its endpoints. This observation allows us to assume without loss of generality that the edge weights on the input graph obey the triangle inequality. Although we have been unable to exploit the triangle inequality property in our algorithms, it is possible that this remark holds the key to improving our results. In this direction, Garg and Hochbaum [17] have recently given an $O(\log k)$-approximation algorithm for the $k$MST problem for points on the plane using an extension of our lower-bounding technique in Section 4.

Table 1 is incomplete. It would be interesting to know the complexity of the minimum-diameter-cost spanning tree problem when the distance values are uniform. Note that any star tree on the nodes provides a 2-approximation to the minimum-diameter-cost spanning tree in this case. The above problem can be shown to be polynomial-time equivalent to the following tree reconstruction problem: given integral nonnegative distances $d_{ij}$ for every pair of vertices $i,j$, does there exist a spanning tree on these nodes such that the distance between $i$ and $j$ in the tree is at most $d_{ij}$?

## 7.2 Maximum acyclic subgraph

In the course of our research we considered the $k$-forest problem: given an undirected graph is there a set of $k$ nodes that induces an acyclic subgraph? The optimization version of this problem is the maximum acyclic subgraph problem. Since this problem is complementary to the minimum feedback vertex set problem [15], NP-completeness follows. While the feedback vertex set problem is 4-approximable [6], we can show that the maximum acyclic subgraph problem is hard to approximate within a reasonable factor using an approximation-preserving transformation from the maximum independent set problem [5]. This same result has also been derived in a more general form in [24].

**Theorem 7.1** *There is a constant $\epsilon > 0$ such that the maximum acyclic subgraph problem cannot be approximated within a factor $\Omega(n^\epsilon)$ unless $P = NP$.*

**Proof:** Note that any acyclic subgraph of size $S$ contains a maximum independent set of size at least $S/2$, since acyclic subgraphs are bipartite and each partition is an independent set. Further, every independent set is also an acyclic subgraph. These two facts show that the existence of a $\rho$-approximation algorithm for the maximum acyclic subgraph problem implies the existence of a $2\rho$-approximation algorithm for the maximum independent set problem. But by the result in [5] we know that there is a constant $\epsilon > 0$ such that the maximum independent set problem cannot be approximated within a factor $\Omega(n^\epsilon)$ unless $P = NP$. Hence, the same is true of the maximum acyclic subgraph problem. □

**Acknowledgements:** The authors wish to thank Alex Zelikovsky and Naveen Garg for helpful conversations during the initial stages of the paper and Professor John Oomen for pointing out that our algorithm for points in the plane extends to the rectilinear case too.